\documentclass{amsart}
\usepackage{amssymb}

\theoremstyle{definition}

\theoremstyle{remark}

\numberwithin{equation}{section}

\input{tcilatex}

\begin{document}
\title{$L_{p}$ regularity for convolution operator equations in Banach spaces%
}
\author{Rishad Shahmurov}
\address{Department of Mathematics, Yeditepe University, Kayishdagi Caddesi,
34755 Kayishdagi, Istanbul, Turkey}
\email{shahmurov@hotmail.com}
\date{}

\begin{abstract}
Here we utilize operator--valued $L_{q}\rightarrow L_{p}$ Fourier multiplier
theorems to establish lower bound estimates for large class of elliptic
integro-differential equations in $R^{d}.$ Moreover, we investigate
separability properties of parabolic convolution operator equations that
arise in heat conduction problems in materials with fading memory. Finally,
we give some remarks on optimal regularity of elliptic differential
equations and Cauchy problem for parabolic equations.
\end{abstract}

\maketitle




\section*{1. Introduction, notations and background}

After pioneering works of Herbert Amann [1] and Lutz Weis [26] on
operator-valued Fourier multiplier theorems (OFMT), theory of
differential--operator equations (DOEs) in Banach valued function spaces is
improved significantly. Many researchers applied them in the investigation
of different classes of equations especially in maximal $L_{p}$ and $%
B_{p,r}^{s}$ regularity for parabolic and elliptic DOE. The exposition of
FMT, their applications and some related references can be found in $\left[ 1%
\right] ,$ $\left[ 3\right] ,$ $\left[ 5\right] ,$ $\left[ 8-11\right] ,$
and $\left[ 26\right] $. For the references concerning FMT in periodic
function spaces, optimal regularity results for convolution operator
equations (COE) and delay DOE see e.g. [12-14], [18], [21] and reference
therein.

Here we shall first extend the well known FMT in [27] and [23]. Then we will
establish lower bound estimates for elliptic type integro-differential
operators of form 
\begin{equation*}
Lu=\ \sum\limits_{k,j=1}^{d}c_{kj}\frac{\partial ^{2}u(x)}{\partial
x_{k}\partial x_{j}}+\sum\limits_{k,j=1}^{d}a_{kj}\ast \frac{\partial
^{2}u(x)}{\partial x_{k}\partial x_{j}}+b_{0}Au+b_{1}\ast Au\eqno(1.1)
\end{equation*}%
where $c_{kj},b_{0}\in \mathbf{C,}$ $a_{kj},b_{1}\in S^{\prime }(R^{d},%
\mathbf{C})$ and $A$ is\ a possible unbounded operator in a Banach space $E$%
. Particularly, we show the following Sobolev type estimates%
\begin{equation*}
\left\Vert u\right\Vert _{L_{p}(R^{d},E)}\leq C\left\Vert Lu\right\Vert
_{L_{q}(R^{d},E)}\eqno(1.2)
\end{equation*}%
for exponents $1<q\leq p<\infty $ satisfying gap condition 
\begin{equation*}
\frac{1}{q}-\frac{1}{p}\leq \frac{2}{d}.
\end{equation*}%
Next we will prove separability for parabolic COE 
\begin{equation*}
Lu=\ a_{0}u^{\prime }+a_{1}\ast u^{\prime }+b_{0}Au+b_{1}\ast Au=f\eqno(1.3)
\end{equation*}%
where $f\in X=L_{p}(R,E),$ $a_{0},b_{0}\in \mathbf{C,}$ $a_{1},b_{1}\in
S^{\prime }(R,\mathbf{C})$ and $A$ is a possible unbounded operator in a
Banach space $E$. Our aim is to obtain the following coercive estimate for
(1.3)%
\begin{equation*}
\left\Vert u^{\prime }\right\Vert _{X}+\left\Vert a_{1}\ast u^{\prime
}\right\Vert _{X}+\left\Vert Au\right\Vert _{X}+\left\Vert b_{1}\ast
Au\right\Vert _{X}\leq C\left\Vert f\right\Vert _{X}.\eqno(1.4)
\end{equation*}%
Note that, model problems for (1.3) are heat conduction problems with fading
memory, population dynamic problems, etc. These problems can be inferred
from (1.3) by choosing $A$ as a second order elliptic differential operator
and $E=L_{r}(R^{2})$ along with appropriate boundary conditions. For
physical interpretations and for detailed information about problems with
fading memory see e.g. [15] and [19].

In the last section we will give some remarks on well-known results
concerning parabolic and elliptic DOE.

Let $\alpha =(\alpha _{1},\alpha _{2},\cdots ,\alpha _{n}),$ where $\alpha
_{i}$ are integers. An $E$--valued generalized function $D^{\alpha }f$ is
called a generalized derivative in the sense of Schwartz distributions, if
the equality 
\begin{equation*}
<D^{\alpha }f,\varphi >=(-1)^{|\alpha |}<f,D^{\alpha }\varphi >
\end{equation*}%
holds for all $\varphi \in S.$

\vspace{3mm}The Fourier transform $F:S(X)\rightarrow S(X)$ is defined by 
\begin{equation*}
(Ff)(t)\equiv \hat{f}(t)=\int\limits_{R^{N}}\exp (-its)f(s)ds
\end{equation*}%
is an isomorphism whose inverse is given by 
\begin{equation*}
(F^{-1}f)(t)\equiv \check{f}(t)=(2\pi )^{-N}\int\limits_{R^{N}}\exp
(its)f(s)ds,
\end{equation*}%
where $f\in S(X)$ and $t\in R^{N}.$ It is clear that 
\begin{equation*}
F(D_{x}^{\alpha }f)=(i\xi _{1})^{\alpha _{1}}\cdots (i\xi _{n})^{\alpha _{n}}%
\hat{f},~D_{\xi }^{\alpha }(F(f))=F[(-ix_{n})^{\alpha _{1}}\cdots
(-ix_{n})^{\alpha _{n}}f]
\end{equation*}%
for all $f\in S^{\dagger }(R^{n};E).$

Let $\mathbf{C}$ be a set of complex numbers and

\begin{equation*}
\ S_{\varphi }=\left\{ \xi ;\text{ }\xi \in \mathbf{C}\text{, \ }\left\vert
\arg \xi \right\vert \leq \varphi \right\} \text{, }0\leq \varphi <\pi .
\end{equation*}%
Suppose $E_{1}$ and\ $E_{2}$ are two Banach spaces. $B\left(
E_{1},E_{2}\right) $ will denote the space of all bounded linear operators
from $E_{1}$ to $E_{2}.$

A linear operator\ $A$ is said to be $\varphi $-positive in\ a Banach\ space 
$E$,\ with bound $M$ if\ $D\left( A\right) $ is dense in $E$ and 
\begin{equation*}
\left\Vert \left( A+\lambda I\right) ^{-1}\right\Vert _{B\left( E\right)
}\leq M\left( 1+\left\vert \lambda \right\vert \right) ^{-1}
\end{equation*}%
for all $\lambda \in S_{\varphi },$ with $\varphi \in \left[ 0,\left. \pi
\right) \right. ,$where $M$ is a positive constant and $I$ is identity
operator in $E.$

$E\left( A^{\theta }\right) $ denotes the space $D\left( A^{\theta }\right) $
with graphical norm 
\begin{equation*}
\left\Vert u\right\Vert _{E\left( A^{\theta }\right) }=\left( \left\Vert
u\right\Vert ^{p}+\left\Vert A^{\theta }u\right\Vert ^{p}\right) ^{\frac{1}{p%
}},1\leq p<\infty ,-\infty <\theta <\infty .
\end{equation*}%
We indicate mixed derivative in the following form 
\begin{equation*}
D^{\alpha }=D_{1}^{\alpha _{1}}D_{2}^{\alpha _{2}}\cdots D_{n}^{\alpha
_{n}},~D_{k}^{i}=\left( \frac{\partial }{\partial x_{k}}\right) ^{i}.
\end{equation*}%
Suppose $\Omega \subset R^{n}.$ Then $W_{p}^{l}\left( \Omega ;E_{0},E\right) 
$ is a space of functions $u\in L_{p}\left( \Omega ;E_{0}\right) $ such that 
$D_{k}^{l}u=\frac{\partial ^{^{l}}u}{\partial x_{k}^{^{l}}}\in L_{p}\left(
\Omega ;E\right) $ and 
\begin{equation*}
\left\Vert u\right\Vert _{W_{p}^{l}\left( \Omega ;E_{0},E\right)
}=\left\Vert u\right\Vert _{L_{p}\left( \Omega ;E_{0}\right)
}+\sum\limits_{k=1}^{n}\left\Vert D_{k}^{l}u\right\Vert _{L_{p}\left( \Omega
;E\right) }<\infty .
\end{equation*}%
For $E_{0}=E$ the space $W_{p}^{l}\left( \Omega ;E_{0},E\right) $ will be
denoted by $W_{p}^{l}\left( \Omega ;E\right) .$

\vspace{3mm}

\section*{2. $L_{q}\rightarrow L_{p}$ FMT}

In this section we shall study scalar-valued FMT from $L_{q}(X)$ to $%
L_{p}(X) $ for $1<q\leq p<\infty $. Let us first introduce some basic
definitions and facts.

\vspace{3mm}

\textbf{Definition 2.0}$.$\textbf{\ }A Banach space $X$ is called UMD space
if $X$-valued martingale difference sequences are unconditional in $%
L_{p}\left( R^{d};X\right) $ for $p\in \left( 1,\infty \right) .$ i.e. there
exists a positive constant $C_{p}$ such that for any martingale $\left\{
f_{k},k\in \mathbf{N}_{0}\right\} $ (see $\left[ \text{6, \S\ 5.}\right] $
[17]), any choice of signs $\left\{ \varepsilon _{k},k\in \mathbf{N}\right\}
\in \left\{ -1,1\right\} $ and $N\in \mathbf{N}$%
\begin{equation*}
\left\Vert f_{0}+\sum\limits_{k=1}^{N}\varepsilon _{k}\left(
f_{k}-f_{k-1}\right) \right\Vert _{L_{p}\left( \Omega ,\Sigma ,\mu ,X\right)
}\leq C_{p}\left\Vert f_{N}\right\Vert _{L_{p}\left( \Omega ,\Sigma ,\mu
,X\right) }.
\end{equation*}%
It is shown in [2] and $\left[ 4\right] $ that\ a Hilbert operator 
\begin{equation*}
\left( Hf\right) \left( x\right) =\lim\limits_{\varepsilon \rightarrow
0}\int\limits_{\left\vert x-y\right\vert >\varepsilon }\frac{f\left(
y\right) }{x-y}dy
\end{equation*}%
is bounded in the space $L_{p}\left( R,X\right) ,$ $p\in \left( 1,\infty
\right) $ for only those spaces $X,$ which possess the UMD property. UMD
spaces include e.g. $L_{p}$, $l_{p}$ spaces and Lorentz spaces $L_{pq},$ $p,$
$q\in \left( 1,\infty \right) $.

\vspace{3mm}

\textbf{Definition 2.1}$.$ Let $X$ and $Y$ be Banach spaces. A family of
operators $\tau \subset B\left( X,Y\right) $ is called $R$-bounded (see e.g. 
$\left[ 5\right] $ and [11]) if there is a positive constant $C$ and $p\in
\lbrack 1,\infty )$ such that for each $N\in \mathbf{N,}$ $T_{j}\in \tau $, $%
x_{j}\in X$ and for all independent, symmetric, $\left\{ -1,1\right\} -$%
valued random variables $r_{j}$ on a probability space $(\Omega ,\Sigma ,\mu
)$ the inequality%
\begin{equation*}
\left\Vert \sum\limits_{j=1}^{N}r_{j}T_{j}x_{j}\right\Vert _{L_{p}(\Omega
,Y)}\leq C\left\Vert \sum\limits_{j=1}^{N}r_{j}x_{j}\right\Vert
_{L_{p}(\Omega ,X)},
\end{equation*}%
is valid. The smallest such $C$ is called $R$-bound of $\tau ,$ we denote it
by $R_{p}(\tau ).$

Let us note that wide classes of classical operators are $R$-bounded. (see
[10] and reference therein). The basic properties of $R$-boundedness are
collected in the recent monograph of Denk et al. [5]. For the reader's
convenience, we present some results from [5].

\vspace{3mm}

(a)The definition of $R$-boundedness is independent of $p\in \lbrack
1,\infty ).$

(b) If $\tau \subset B\left( X,Y\right) $ is $R$-bounded then it is
uniformly bounded with%
\begin{equation*}
\sup \left\{ \left\Vert T\right\Vert :\text{ }T\in \tau \right\} \leq
R_{p}(\tau ).
\end{equation*}

(c) If $X$ and $Y$ are Hilbert spaces, $\tau \subset B\left( X,Y\right) $ is 
$R$-bounded $\Longleftrightarrow \tau $ is uniformly bounded.

(d) Let $X,$ $Y$ be Banach spaces and $\tau _{1}$, $\tau _{2}\subset B(X,Y)$
be $R$-bounded. Then%
\begin{equation*}
\tau _{1}+\tau _{2}=\left\{ T+S:T\in \tau _{1},\text{ }S\in \tau _{2}\right\}
\end{equation*}%
is $R$-bounded as well, and $R_{p}(\tau _{1}+\tau _{2})\leq R_{p}(\tau
_{1})+R_{p}(\tau _{2})$ .

(e) Let $X,$ $Y,$ $Z$ be Banach spaces and $\tau _{1}\subset B(X,Y)$ and $%
\tau _{2}\subset B(Y,Z)$ be $R$-bounded. Then%
\begin{equation*}
\tau _{1}\tau _{2}=\left\{ ST:T\in \tau _{1},\text{ }S\in \tau _{2}\right\}
\end{equation*}%
is $R$-bounded as well, and $R_{p}(\tau _{1}\tau _{2})\leq R_{p}(\tau
_{1})R_{p}(\tau _{2})$.

One of the most important tools in $R$-boundedness is the contraction
principle of Kahane. We shall frequently apply it in the next sections.

\vspace{3mm}

\textbf{[5}, \textbf{Lemma 3.5.] }Let $X$ be a Banach spaces, $n\in N,$ $%
x_{j}\in X,$ $r_{j}$ independent, symmetric, $\left\{ -1,1\right\} $-valued
random variables on a probability space $(\Omega ,\Sigma ,\mu )$ and $\alpha
_{j,}\beta _{j}\in \mathbf{C}$ such that $\left\vert \alpha _{j}\right\vert
\leq \left\vert \beta _{j}\right\vert ,$ for each $j=1,\cdot \cdot \cdot ,N.$
Then%
\begin{equation*}
\left\Vert \sum\limits_{j=1}^{N}\alpha _{j}r_{j}x_{j}\right\Vert
_{L_{p}(\Omega ,X)}\leq 2\left\Vert \sum\limits_{j=1}^{N}\beta
_{j}r_{j}x_{j}\right\Vert _{L_{p}(\Omega ,X)}.
\end{equation*}%
The constant 2 can be omitted in case where $\alpha _{j}$ and $\beta _{j}$
are real.

\vspace{3mm}

\textbf{Theorem 2.2. }Let $X$ be an UMD space and $1<q\leq p<\infty .$ If
for a bounded function $\psi :R^{d}\backslash \left\{ 0\right\} \rightarrow 
\mathbf{C}$ 
\begin{equation*}
\sup \left\{ \left\vert \xi \right\vert ^{\left\vert \alpha \right\vert
+d\left( \frac{1}{q}-\frac{1}{p}\right) }\left\vert D^{\alpha }\psi (\xi
)\right\vert ,\text{ }\xi \in R^{d}\backslash \left\{ 0\right\} ,\text{ }%
\alpha \leq (1,...,1)\right\} <\infty \text{ }\eqno(2.1)
\end{equation*}%
then 
\begin{equation*}
\left\Vert T_{\psi }f\right\Vert _{L_{p}(R^{d},Y)}=\left\Vert F^{-1}\left[
\psi (\cdot )\hat{f}(\cdot )\right] \right\Vert _{L_{p}(R^{d},Y)}\leq
C\left\Vert f\right\Vert _{L_{q}(R^{d},X)}
\end{equation*}%
for all $f\in S(R^{d},X).$

\textbf{Proof. }The main idea is to apply $L_{p}\rightarrow L_{p}$ FMT along
with Sobolev embedding and to use nice properties of a function $%
s(t)=\left\vert t\right\vert ^{d(\frac{1}{q}-\frac{1}{p})}$. Since $\psi $
satisfies the general Miklin's condition i.e.%
\begin{equation*}
\sup \left\{ \left\vert \xi \right\vert ^{\left\vert \alpha \right\vert
+d\left( \frac{1}{q}-\frac{1}{p}\right) }\left\vert D^{\alpha }\psi (\xi
)\right\vert ,\text{ }\xi \in R^{d}\backslash \left\{ 0\right\} ,\text{ }%
\alpha \leq (1,...,1)\right\} \leq C,
\end{equation*}%
$s\cdot \psi $ satisfies classical one i.e. 
\begin{equation*}
\sup \left\{ \left\vert \xi \right\vert ^{\left\vert \alpha \right\vert
}\left\vert D^{\alpha }(s\psi )(\xi )\right\vert ,\text{ }\xi \in
R^{d}\backslash \left\{ 0\right\} ,\text{ }\alpha \leq (1,...,1)\right\}
\leq C.
\end{equation*}%
Therefore applying the Sobolev embedding theorem and by using the definition
of homogeneous Sobolev spaces (in the sense of Riesz potentials) we get
desired result: 
\begin{equation*}
\begin{array}{lll}
\left\Vert F^{-1}\left[ \psi \hat{f}\right] \right\Vert _{L_{p}\left(
R^{d};X\right) } & \leq & \displaystyle K\left\Vert F^{-1}\left[ \psi \hat{f}%
\right] \right\Vert _{\dot{W}_{q}^{d\left( \frac{1}{q}-\frac{1}{p}\right)
}\left( R^{d};X\right) } \\ 
&  &  \\ 
& = & \displaystyle K\left\Vert F^{-1}\left[ \left\vert \cdot \right\vert
^{d\left( \frac{1}{q}-\frac{1}{p}\right) }\psi (\cdot )\hat{f}(\cdot )\right]
\right\Vert _{L_{q}\left( R^{d};X\right) } \\ 
&  &  \\ 
& = & \displaystyle K\left\Vert F^{-1}\left[ \left( s\psi \right) \hat{f}%
(\cdot )\right] \right\Vert _{L_{q}\left( R^{d};X\right) }\leq C\left\Vert
f\right\Vert _{L_{q}\left( R^{d};X\right) }.%
\end{array}%
\end{equation*}

\hbox{\vrule height7pt width5pt}

The following result is extended version of operator valued Miklin theorem
in [23]. Theorem 2.3 can be proven in a similar manner as Theorem 2.2.

\vspace{3mm}

\textbf{Theorem 2.3. }Let $X$ and $Y$ be UMD spaces and $1<q\leq p<\infty .$
If for a bounded function $M:R^{d}\backslash \left\{ 0\right\} \rightarrow
B(X,Y),$ 
\begin{equation*}
R\left\{ \left\vert x\right\vert ^{\left\vert \alpha \right\vert +d\left( 
\frac{1}{q}-\frac{1}{p}\right) }D^{\alpha }M(x),\text{ }x\in R^{d}\backslash
\left\{ 0\right\} ,\text{ }(\alpha \leq 1,...,1)\text{ }\right\} <\infty 
\text{ }\eqno(2.2)
\end{equation*}%
then 
\begin{equation*}
\left\Vert T_{M}f\right\Vert _{L_{p}(R^{d},Y)}=\left\Vert F^{-1}\left[
M(\cdot )\hat{f}(\cdot )\right] \right\Vert _{L_{p}(R^{d},Y)}\leq
C\left\Vert f\right\Vert _{L_{q}(R^{d},X)}
\end{equation*}%
for all $f\in S(R^{d},X).$

\section*{3. Sobolev type estimates for (1.1)}

Let us consider the second order elliptic integro-differential equation
(1.1) in $R^{d}$. Here we characterize conditions on coefficients of (1.1)
so that they imply Sobolev type estimate (1.2).

Since we utilize Fourier integral methods we naturally impose conditions on
symbols. Therefore, to avoid contradictions due to Riemann-Lebesgue lemma we
make some auxilary assumptions along with certain regularity and ellipticity
conditions on coefficients of $L$.

\vspace{3mm}

\textbf{Definition 3.0. }Let $E$ be a Banach space and $D\left( A\right) $
dense in $E.$ A $\varphi $-positive operator $A$ is said to be $R$-positive
if the following set%
\begin{equation*}
\left\{ (1+\xi )(A+\xi )^{-1}:\text{ }\xi \in S_{\varphi },\text{ }\varphi
\in \lbrack 0,\pi )\right\}
\end{equation*}%
is $R$-bounded.

In what follows $\dot{R}$ will denote the set of real numbers excluding zero
i.e. $\dot{R}^{d}=R^{d}/\left\{ 0\right\} .$ 

\vspace{3mm}

\textbf{Condition 3.1. }Suppose the following are satisfied:

(1) $c_{kj},b_{0}\in \mathbf{C,}$ $a_{kj},b_{1}\in S^{\prime }(R^{d},\mathbf{%
C}),$ $\hat{a}_{kj},\hat{b}_{1}\in C^{d}(\dot{R}^{d},\mathbf{C})$ and 
\begin{equation*}
C_{b}=\inf_{\xi \in \dot{R}^{d}}\left\vert b_{0}+\hat{b}_{1}(\xi
)\right\vert >0\text{ };
\end{equation*}

(2) There exists a constant $C$ such that%
\begin{equation*}
\left\vert N\left( \xi \right) \right\vert =\left\vert
\sum\limits_{k,j=1}^{d}\left( c_{kj}+\hat{a}_{kj}(\xi )\right) \xi _{k}\xi
_{j}\right\vert \geq C\left\vert \xi \right\vert ^{2};
\end{equation*}

(3)%
\begin{equation*}
\eta \left( \xi \right) =\frac{N\left( \xi \right) }{\hat{b}_{1}(\xi )+b_{0}}%
\in S_{\varphi }\text{ for }\varphi \in \left[ 0,\right. \left. \pi \right) ;
\end{equation*}

(4) there are some constants $C_{i}$ such that for all $\xi \in \dot{R}^{d}$ 
\begin{eqnarray*}
\left\vert \xi \right\vert ^{\left\vert \beta \right\vert }\left\vert \frac{%
\partial ^{\left\vert \beta \right\vert }}{\partial \xi _{1}^{\beta
_{1}}\partial \xi _{2}^{\beta 2}\cdot \cdot \cdot \partial \xi _{d}^{\beta
_{d}}}\hat{a}_{kj}(\xi )\right\vert  &\leq &C_{0}\text{ for all }k,\text{ }%
j=1,2,\cdot \cdot \cdot ,d\text{,} \\
&& \\
\left\vert \xi \right\vert ^{\left\vert \beta \right\vert }\left\vert \frac{%
\partial ^{\left\vert \beta \right\vert }}{\partial \xi _{1}^{\beta
_{1}}\partial \xi _{2}^{\beta 2}\cdot \cdot \cdot \partial \xi _{d}^{\beta
_{d}}}\hat{b}_{1}(\xi )\right\vert  &\leq &C_{1},
\end{eqnarray*}%
where $\beta _{i}\in \left\{ 0,1\right\} $ and $0\leq \left\vert \beta
\right\vert \leq d.$ \vspace{3mm}

\textbf{Theorem 3.2.}\ Suppose $E$ is an UMD space and Condition 3.1 holds.
Let $A$ be an $R$-positive operator in $E$ with $0\leq \varphi <\pi $. Then, 
$\left( 1.1\right) $ satisfies Sobolev type estimate (1.2) for exponents $%
1<q\leq p<\infty $ satisfying gap condition 
\begin{equation*}
\frac{1}{q}-\frac{1}{p}\leq \frac{2}{d}.
\end{equation*}

To prove our main result we will need the following preliminary lemmas.

\vspace{3mm}

\textbf{Lemma 3.3. }Let $A$ be an $R$-positive operator in $E$ and assume
Condition 3.1 holds. Then the following set 
\begin{equation*}
\left\{ \left\vert \xi \right\vert ^{d\left( \frac{1}{q}-\frac{1}{p}\right)
}\sigma (\xi );\text{ }\xi \in \dot{R}^{d}\right\} 
\end{equation*}%
is $R$-bounded where%
\begin{equation*}
\sigma (\xi )=\frac{1}{\hat{b}_{1}(\xi )+b_{0}}\left( A+\eta \left( \xi
\right) \right) ^{-1}.
\end{equation*}

\textbf{Proof.} Since $\eta \left( \xi \right) \in S_{\varphi }$ for $%
\varphi \in \left[ 0,\right. \left. \pi \right) ,$ from [7, Lemma 2.3] there
exist $K>0$ independent of $\xi $ so that 
\begin{equation*}
\left\vert 1+\eta \left( \xi \right) \right\vert ^{-1}\leq K(1+\left\vert
\eta \left( \xi \right) \right\vert )^{-1}.
\end{equation*}%
Therefore, for all $\xi \in \dot{R}^{d}$ we have uniform estimate 
\begin{equation*}
\frac{\left\vert \xi \right\vert ^{d\left( \frac{1}{q}-\frac{1}{p}\right) }}{%
\left\vert \hat{b}_{1}(\xi )+b_{0}\right\vert }\times \frac{1}{\left\vert
1+\eta \left( \xi \right) \right\vert }\leq K\frac{\left\vert \xi
\right\vert ^{d\left( \frac{1}{q}-\frac{1}{p}\right) }}{\left\vert \hat{b}%
_{1}(\xi _{j})+b_{0}\right\vert +\left\vert N\left( \xi \right) \right\vert }%
\leq K\frac{\left\vert \xi \right\vert ^{d\left( \frac{1}{q}-\frac{1}{p}%
\right) }}{C_{b}+C\left\vert \xi \right\vert ^{2}}\leq \frac{K}{C}.
\end{equation*}%
Now let us define families of operators 
\begin{equation*}
\tau _{1}=\left\{ T_{j}=\left\vert \xi ^{j}\right\vert ^{d\left( \frac{1}{q}-%
\frac{1}{p}\right) }\sigma (\xi ^{j});\text{ }\xi ^{j}\in \dot{R}%
^{d}\right\} 
\end{equation*}%
and 
\begin{equation*}
\tau _{2}=\left\{ S_{j}=(1+\eta \left( \xi ^{j}\right) )\left( A+\eta \left(
\xi ^{j}\right) \right) ^{-1};\text{ }\xi ^{j}\in \dot{R}^{d}\right\} .
\end{equation*}%
Taking into consideration $R$-positivity of $A$, applying assumptions of
Condition 3.1 and Kahane's contraction principle [7, Lemma 3.5] we get
desired result:%
\begin{equation*}
\begin{array}{lll}
\left\Vert \sum\limits_{j=1}^{N}r_{j}T_{j}x_{j}\right\Vert _{X} & = & %
\displaystyle\left\Vert \sum\limits_{j=1}^{N}r_{j}\frac{\left\vert \xi
^{j}\right\vert ^{d\left( \frac{1}{q}-\frac{1}{p}\right) }}{(\hat{b}_{1}(\xi
^{j})+b_{0})(1+\eta \left( \xi ^{j}\right) )}S_{j}x_{j}\right\Vert _{X} \\ 
&  &  \\ 
& \leq  & \displaystyle2\frac{K}{C}\left\Vert
\sum\limits_{j=1}^{N}r_{j}S_{j}x_{j}\right\Vert _{X}\leq 2\frac{K}{C}%
R_{p}(\tau _{2})\left\Vert \sum\limits_{j=1}^{N}r_{j}x_{j}\right\Vert _{X}.%
\end{array}%
\end{equation*}%
where $X=L_{p}((0,1),E).$ Hence 
\begin{equation*}
R_{p}(\tau _{1})\leq 2KR_{p}(\tau _{2}).
\end{equation*}%
\hbox{\vrule
height7pt width5pt}

In the next lemmas we will estimate $R$-bounds of partial derivatives of $%
\sigma (\xi ).$

\vspace{3mm}

\textbf{Lemma 3.4. }Let $A$ be an $R$-positive operator in $E$ and assume
Condition 3.1 holds. Then, the following set 
\begin{equation*}
\left\{ \text{ }\left\vert \xi \right\vert ^{1+d\left( \frac{1}{q}-\frac{1}{p%
}\right) }\frac{\partial }{\partial \xi _{i}}\sigma (\xi ),\text{ }\xi \in 
\dot{R}^{d}\text{ }\right\} 
\end{equation*}%
is $R$-bounded.

\textbf{Proof. }It clear to see that first derivative of $\sigma (\xi )$
consist of 3 terms namely%
\begin{equation*}
\begin{array}{lll}
\sigma _{1}(\xi ) & = & \displaystyle\frac{-\frac{\partial \hat{b}_{1}}{%
\partial \xi _{i}}(\xi )}{\left( \hat{b}_{1}(\xi )+b_{0}\right) ^{2}}\left(
A+\eta \left( \xi \right) \right) ^{-1}, \\ 
&  &  \\ 
\sigma _{2}(\xi ) & = & \displaystyle\left( \frac{1}{\hat{b}_{1}(\xi )+b_{0}}%
\right) ^{3}\frac{\partial N\left( \xi \right) }{\partial \xi _{i}}\hat{b}%
_{1}(\xi )N\left( \xi \right) \left( A+\eta \left( \xi \right) \right) ^{-2},
\\ 
&  &  \\ 
\sigma _{3}(\xi ) & = & \displaystyle\left( \frac{1}{\hat{b}_{1}(\xi )+b_{0}}%
\right) ^{3}\frac{\partial \hat{b}_{1}(\xi )}{\partial \xi _{i}}N\left( \xi
\right) \left( A+\eta \left( \xi \right) \right) ^{-2}.%
\end{array}%
\end{equation*}%
For the sake of simplicity we will only estimate $R$-bound of the set 
\begin{equation*}
\left\{ \left\vert \xi \right\vert ^{1+d\left( \frac{1}{q}-\frac{1}{p}%
\right) }\sigma _{1}(\xi ):\xi \in \dot{R}^{d}\right\} .
\end{equation*}%
Define a family of operators 
\begin{equation*}
\tau =\left\{ T_{j}=\left\vert \xi ^{j}\right\vert ^{1+d\left( \frac{1}{q}-%
\frac{1}{p}\right) }\sigma _{1}(\xi ^{j});\text{ }\xi ^{j}\in \dot{R}%
^{d}\right\} .
\end{equation*}%
Making use of Condition 3.1 and $R$-positivity of $A$ we get%
\begin{equation*}
\begin{array}{lll}
\left\Vert \sum\limits_{j=1}^{N}r_{j}T_{j}x_{j}\right\Vert _{X} & = & %
\displaystyle\left\Vert \sum\limits_{j=1}^{N}r_{j}\frac{-\left\vert \xi
^{j}\right\vert \frac{\partial \hat{b}_{1}}{\partial \xi _{i}^{j}}(\xi )}{%
\hat{b}_{1}(\xi ^{j})+b_{0}}\frac{\left\vert \xi ^{j}\right\vert ^{d\left( 
\frac{1}{q}-\frac{1}{p}\right) }}{\left( \hat{b}_{1}(\xi ^{j})+b_{0}\right)
(1+\eta \left( \xi ^{j}\right) )}S_{j}x_{j}\right\Vert _{X} \\ 
&  &  \\ 
& \leq  & \displaystyle2\frac{K}{C}\frac{C_{1}}{C_{b}}\left\Vert
\sum\limits_{j=1}^{N}r_{j}S_{j}x_{j}\right\Vert _{X}\leq 2\frac{K}{C}\frac{%
C_{1}}{C_{b}}R_{p}(\tau _{2})\left\Vert
\sum\limits_{j=1}^{N}r_{j}x_{j}\right\Vert _{X},%
\end{array}%
\end{equation*}%
which implies 
\begin{equation*}
R_{p}(\tau )\leq 2\frac{K}{C}\frac{C_{1}}{C_{b}}R_{p}(\tau _{2}).
\end{equation*}%
In a similar fashion one can also prove%
\begin{equation*}
R\left\{ \left\vert \xi \right\vert ^{1+d\left( \frac{1}{q}-\frac{1}{p}%
\right) }\sigma _{j}(\xi ):\xi \in \dot{R}^{d}\right\} \leq M.
\end{equation*}%
Hence we get assertion of the lemma. \ \hbox{\vrule
height7pt width5pt}

The next result is generalization of Lemma 3.4. We will omit the prove since
it analogously follows from previous results.

\vspace{3mm}

\textbf{Lemma 3.5. }Assume $A$ is an $R$-positive operator in $E$ and
Condition 3.1 holds. Then, the following set 
\begin{equation*}
\left\{ \text{ }\left\vert \xi \right\vert ^{\left\vert \alpha \right\vert
+d\left( \frac{1}{q}-\frac{1}{p}\right) }D^{\alpha }\sigma (\xi ),\text{ }%
\xi \in \dot{R}^{d},\text{ }\alpha \leq ^{\text{comp.wise}}(1,...,1)\text{ }%
\right\} 
\end{equation*}%
is $R$-bounded.

\textbf{Proof of Theorem 3.2}. Suppose $Lu=f$ for some $f\in L_{q}\left(
R^{d};E\right) .$ Taking into consideration Condition 3.1 and applying the
Fourier transform to both side of equation\ we get 
\begin{equation*}
u\left( x\right) =F^{-1}\left[ \frac{1}{\hat{b}_{1}(\xi )+b_{0}}\left( A+%
\frac{N\left( \xi \right) }{\hat{b}_{1}(\xi )+b_{0}}\right) ^{-1}\hat{f}%
\left( \xi \right) \right] .
\end{equation*}%
Now it is easy to see that (1.2) is equivalent to $L_{q}\left(
R^{d};E\right) \rightarrow L_{p}\left( R^{d};E\right) $ boundedness of above
Fourier multiplier operator. Since the operator-valued multiplier function%
\begin{equation*}
m(\xi )=\frac{1}{\hat{b}_{1}(\xi )+b_{0}}\left( A+\frac{N\left( \xi \right) 
}{\hat{b}_{1}(\xi )+b_{0}}\right) ^{-1}
\end{equation*}%
satisfies assumptions of Theorem 2.3 we complete the proof. \ \ 
\hbox{\vrule
height7pt width5pt}

\section*{4. Optimal regular Parabolic COE}

Due to its nice applications, many researchers investigated (1.3) in various
function spaces. The Besov space regularity for (1.3) is studied in [1] and
the Holder space ($C^{\alpha }$ with $0<\alpha <1)$ case is\ presented in
[14]$.$ Moreover, maximal regularity results for (1.3) in different periodic
function spaces can be found in very recent paper [12].

Here we study the same problem in $L_{p}(R,E)$ under some natural
assumptions on coefficients$.$ The main tool we implement here will be FMT
of Weis [26]. First we state our assumptions:

\vspace{3mm}

\textbf{Condition 4.1. }Suppose the following are satisfied:

(1) $a_{0},b_{0}\in \mathbf{C,}$ $a_{1},b_{1}\in S^{\prime }(R,\mathbf{C}),$ 
$\hat{a}_{1},\hat{b}_{1}\in C^{1}(\dot{R},\mathbf{C})$ and 
\begin{equation*}
\lim_{\left\vert \xi \right\vert \rightarrow \infty }\inf \left\vert a_{0}+%
\hat{a}_{1}(\xi )\right\vert >0\text{ and }\inf_{\xi \in \dot{R}}\left\vert
b_{0}+\hat{b}_{1}(\xi )\right\vert =C_{b}>0\text{ };
\end{equation*}

(2)%
\begin{equation*}
\frac{i\xi \left( \hat{a}_{1}(\xi )+a_{0}\right) }{\hat{b}_{1}(\xi )+b_{0}}%
\in S_{\varphi }\text{ for }\varphi \in \left[ 0,\right. \left. \pi \right) ;
\end{equation*}

(3) There are constants $C_{i}$ such that 
\begin{eqnarray*}
\left\vert \hat{a}_{1}(\xi )\right\vert  &\leq &C_{0}\text{, }\left\vert \xi 
\frac{d}{d\xi }\hat{a}_{1}(\xi )\right\vert \leq C_{1} \\
&& \\
\left\vert \hat{b}_{1}(\xi )\right\vert  &\leq &C_{2},\text{ }\left\vert \xi 
\frac{d}{d\xi }\hat{b}_{1}(\xi )\right\vert \leq C_{3},\text{ for all }\xi
\in \dot{R}.
\end{eqnarray*}

\vspace{3mm}

\textbf{Theorem 4.2.}\ Assume $E$ is an UMD space and Condition 4.1 holds.
Let $A$ be an $R$-positive operator in $E$ with $0\leq \varphi <\pi $ and $%
1<p<\infty $. Then, the equation $\left( 1.3\right) $ has a unique solution $%
u\in W_{p}^{1}(R,E(A);E)$ satisfying (1.4).

\ \textbf{Proof}. Taking into consideration Condition 4.1 and applying the
Fourier transform to both side of (1.3)\ we get 
\begin{equation*}
u\left( x\right) =F^{-1}\left[ \mu (\xi )\left( A+\eta (\xi )\right) ^{-1}%
\hat{f}\left( \xi \right) \right] .
\end{equation*}%
where%
\begin{equation*}
\eta (\xi )=i\xi \left( \hat{a}_{1}(\xi )+a_{0}\right) \mu (\xi )
\end{equation*}%
and%
\begin{equation*}
\mu (\xi )=\frac{1}{\hat{b}_{1}(\xi )+b_{0}}.
\end{equation*}%
Since%
\begin{equation*}
\left\Vert u\right\Vert _{L_{p}(R,E)}\equiv \left\Vert F^{-1}\left[ \mu (\xi
)\left( A+\eta (\xi )\right) ^{-1}\hat{f}\left( \xi \right) \right]
\right\Vert _{L_{p}(R,E)},
\end{equation*}%
\begin{eqnarray*}
\left\Vert u^{\prime }\right\Vert _{L_{p}(R,E)} &\equiv &\left\Vert F^{-1} 
\left[ i\xi \mu (\xi )\left( A+\eta (\xi )\right) ^{-1}\hat{f}\left( \xi
\right) \right] \right\Vert _{L_{p}(R,E)},\text{ } \\
\left\Vert a_{1}\ast u^{\prime }\right\Vert _{L_{p}(R,E)} &\equiv
&\left\Vert F^{-1}\left[ i\xi \hat{a}_{1}(\xi )\mu (\xi )\left( A+\eta (\xi
)\right) ^{-1}\hat{f}\left( \xi \right) \right] \right\Vert _{L_{p}(R,E)},
\end{eqnarray*}%
\begin{equation*}
\left\Vert Au\right\Vert _{L_{p}(R,E)}\equiv \left\Vert F^{-1}\left[ \mu
(\xi )A\left( A+\eta (\xi )\right) ^{-1}\hat{f}\left( \xi \right) \right]
\right\Vert _{L_{p}(R,E)},
\end{equation*}%
and%
\begin{equation*}
\left\Vert b_{1}\ast Au\right\Vert _{L_{p}(R,E)}\equiv \left\Vert F^{-1} 
\left[ \hat{b}_{1}(\xi )\mu (\xi )A\left( A+\eta (\xi )\right) ^{-1}\hat{f}%
\left( \xi \right) \right] \right\Vert _{L_{p}(R,E)}
\end{equation*}%
it suffices to show 
\begin{equation*}
m_{0}(\xi )=\mu (\xi )\left( A+\eta (\xi )\right) ^{-1},\text{ }m_{1}(\xi
)=i\xi \mu (\xi )\left( A+\eta (\xi )\right) ^{-1}
\end{equation*}%
\begin{equation*}
m_{2}(\xi )=i\xi \hat{a}_{1}(\xi )\mu (\xi )\left( A+\eta (\xi )\right)
^{-1},\text{ }m_{3}(\xi )=\mu (\xi )A\left( A+\eta (\xi )\right) ^{-1}
\end{equation*}%
and%
\begin{equation*}
m_{4}(\xi )=\hat{b}_{1}(\xi )\mu (\xi )A\left( A+\eta (\xi )\right) ^{-1}
\end{equation*}%
are Fourier multipliers in $L_{p}(R,E)$. Therefore, we will prove in several
steps that $m_{i}(\xi )$ satisfy (2.2) for $p=q$ and $d=1.$

\bigskip First we show that $S_{i}=\left\{ m_{i}(\xi ):\xi \in \dot{R}%
\right\} $ are $R$-bounded sets.

\textbf{Lemma 4.3. }Let $E$ be an UMD space and $A$ be an $R$-positive
operator in $E$ with $0\leq \varphi <\pi $. If Condition 4.1 holds then $%
S_{i}=\left\{ m_{i}(\xi ):\xi \in \dot{R}\right\} $ are $R$-bounded sets.

\textbf{Proof. }As in Lemma 3.3 for all $\xi \in \dot{R}$ we have uniform
estimates%
\begin{equation*}
\begin{array}{lll}
\frac{\left\vert \mu (\xi )\right\vert }{\left\vert 1+\eta (\xi )\right\vert 
} & \leq  & \displaystyle\frac{1}{\left\vert \hat{b}_{1}(\xi
)+b_{0}\right\vert +\left\vert i\xi \left( \hat{a}_{1}(\xi )+a_{0}\right)
\right\vert }\leq \frac{1}{C_{b}}, \\ 
&  &  \\ 
\frac{\left\vert i\xi \right\vert \left\vert \mu (\xi )\right\vert }{%
\left\vert 1+\eta (\xi )\right\vert } & \leq  & \displaystyle\frac{%
\left\vert i\xi \right\vert }{\left\vert \hat{a}_{1}(\xi )+a_{0}\right\vert
\left\vert i\xi \right\vert +C_{b}}\leq K, \\ 
&  &  \\ 
\frac{\left\vert i\xi \right\vert \left\vert \mu (\xi )\right\vert
\left\vert \hat{a}_{1}(\xi )\right\vert }{\left\vert 1+\eta (\xi
)\right\vert } & \leq  & \displaystyle KC_{0}, \\ 
&  &  \\ 
\left\vert \mu (\xi )\right\vert \frac{\left\vert \eta (\xi )\right\vert }{%
\left\vert 1+\eta (\xi )\right\vert } & \leq  & \displaystyle\frac{1}{C_{b}}%
\end{array}%
\end{equation*}%
and%
\begin{equation*}
\begin{array}{lll}
\left\vert \hat{b}_{1}(\xi )\right\vert \left\vert \mu (\xi )\right\vert 
\frac{\left\vert \eta (\xi )\right\vert }{\left\vert 1+\eta (\xi
)\right\vert } & \leq  & \displaystyle\frac{C_{2}}{C_{b}}.%
\end{array}%
\end{equation*}%
\textbf{\ } Now let us define families of operators 
\begin{equation*}
\tau _{i}=\left\{ T_{j}^{i}=m_{i}(\xi ^{j});\text{ }\xi ^{j}\in \dot{R}%
\right\} \text{ for }i=0,\cdot \cdot \cdot ,4
\end{equation*}%
and 
\begin{equation*}
\tau =\left\{ S_{j}=(1+\eta \left( \xi ^{j}\right) )\left( A+\eta \left( \xi
^{j}\right) \right) ^{-1};\text{ }\xi ^{j}\in \dot{R}\right\} .
\end{equation*}%
Taking into consideration $R$-positivity of $A$, applying assumptions of
Condition 3.1 and Kahane's contraction principle [7, Lemma 3.5] we get 
\begin{equation*}
\begin{array}{lll}
\left\Vert \sum\limits_{j=1}^{N}r_{j}T_{j}^{0}x_{j}\right\Vert _{X} & \leq 
& \displaystyle\frac{2}{C_{b}}R_{p}(\tau )\left\Vert
\sum\limits_{j=1}^{N}r_{j}x_{j}\right\Vert _{X}, \\ 
&  &  \\ 
\left\Vert \sum\limits_{j=1}^{N}r_{j}T_{j}^{1}x_{j}\right\Vert _{X} & \leq 
& \displaystyle2KR_{p}(\tau )\left\Vert
\sum\limits_{j=1}^{N}r_{j}x_{j}\right\Vert _{X}%
\end{array}%
\end{equation*}%
and%
\begin{equation*}
\begin{array}{lll}
\left\Vert \sum\limits_{j=1}^{N}r_{j}T_{j}^{2}x_{j}\right\Vert _{X} & \leq 
& \displaystyle2KC_{0}R_{p}(\tau )\left\Vert
\sum\limits_{j=1}^{N}r_{j}x_{j}\right\Vert _{X}%
\end{array}%
\end{equation*}%
which implies%
\begin{equation*}
R_{p}(\tau _{0})\leq \frac{2}{C_{b}}R_{p}(\tau ),\text{ }R_{p}(\tau
_{1})\leq 2KR_{p}(\tau )\text{ and }R_{p}(\tau _{2})\leq 2KC_{0}R_{p}(\tau ).
\end{equation*}%
Finally, in view of resolvent properties of positive operators and again by
[7, Lemma 3.5] we deduce%
\begin{equation*}
R_{p}(\tau _{3})=R\left\{ \mu (\xi )\left[ I-\eta (\xi )\left( A+\eta (\xi
)\right) ^{-1}\right] :\text{ }\xi \in \dot{R}\right\} 
\end{equation*}%
\begin{equation*}
\begin{array}{lll}
& \leq  & \displaystyle R\left\{ \mu (\xi )I:\text{ }\xi \in \dot{R}\right\}
+R\left\{ \mu (\xi )\eta (\xi )\left( A+\eta (\xi )\right) ^{-1}\right\}  \\ 
&  &  \\ 
& \leq  & \displaystyle\frac{1}{C_{b}}+\frac{2}{C_{b}}R_{p}(\tau )%
\end{array}%
\end{equation*}%
and 
\begin{equation*}
\begin{array}{lll}
R_{p}(\tau _{4}) & \leq  & \displaystyle\frac{C_{2}}{C_{b}}+\frac{2C_{2}}{%
C_{b}}R_{p}(\tau ).%
\end{array}%
\end{equation*}%
\hbox{\vrule
height7pt width5pt}

Next we will estimate derivatives of operator valued functions $m_{i}(\xi ).$

\vspace{3mm}

\textbf{Lemma 4.4. }Let $E$ be an UMD space and $A$ be an $R$-positive
operator in $E$ with $0\leq \varphi <\pi $. If Condition 4.1 holds then $%
S_{i}=\left\{ \xi \frac{d}{d\xi }m_{i}(\xi ):\xi \in \dot{R}\right\} $ are $R
$-bounded sets.

\textbf{Proof. \ }For the sake of simplicity we shall prove only for $S_{0}.$
The other cases can be proved analogously with the help of above techniques.
Taking derivative of $m_{0}$ and applying similar ideas as in Lemma 4.3 we
get desired result: 
\begin{equation*}
R\left\{ \xi \frac{d}{d\xi }m_{0}(\xi )\right\} \leq R\left\{ \frac{-\xi 
\frac{d}{d\xi }\hat{b}_{1}(\xi )}{\hat{b}_{1}(\xi )+b_{0}}\mu (\xi )\left(
A+\eta (\xi )\right) ^{-1}\right\} +R\left\{ \eta ^{\prime }(\xi )\xi \mu
(\xi )\left( A+\eta (\xi )\right) ^{-2}\right\} 
\end{equation*}%
\begin{equation*}
\leq \frac{2}{C_{b}^{2}}R_{p}(\tau )\left[ C_{3}+K\left( C_{0}+\left\vert
a_{0}\right\vert +\frac{C_{3}\left( C_{0}+\left\vert a_{0}\right\vert
\right) }{C_{b}}+C_{1}\right) \right] <\infty .
\end{equation*}

\hbox{\vrule
height7pt width5pt}

\vspace{3mm}

\textbf{Corollary 4.5. }Let $E$ be an UMD space and $A$ be an $R$-positive
operator in $E$ with $0\leq \varphi <\pi $. If Condition 4.1 holds then $%
m_{i}(\xi )$ are Fourier multipliers in $L_{p}(R,E)$ for $1<p<\infty .$

Since $m_{i}(\xi )$ are Fourier multipliers we complete the proof of Theorem
4.2. \ \ \hbox{\vrule
height7pt width5pt}

\vspace{3mm}

\textbf{Example 4.6. }As an application of our main result we can give a
heat conduction problem in materials with fading memory. Really, choosing $%
E=L_{q}(\Omega ),$ $A=-\frac{d^{2}}{dx^{2}}+c,$ $a_{1}(t)=e^{-m\left\vert
t\right\vert }$, $b_{1}(t)=$ $e^{-k\left\vert t\right\vert }$ in (1.3) we
obtain the following integro-differential equation%
\begin{equation*}
\partial _{t}u+\int\limits_{-\infty }^{\infty }e^{-m\left\vert
t-s\right\vert }\partial _{t}u(s,x)ds=
\end{equation*}%
\begin{equation*}
=f(t,x)+(\partial _{xx}-c)u+\int\limits_{-\infty }^{\infty }e^{-k\left\vert
t-s\right\vert }(\partial _{xx}-c)u(s,x)ds,
\end{equation*}%
\begin{equation*}
u(t,x)_{x\in \partial \Omega }=0,
\end{equation*}%
where $f\in X=L_{p}(R;L_{q}(\Omega ))=L_{p,q}(R\times \Omega )$, $c,m,k>0$
and $\partial \Omega $ is a sufficiently smooth boundary$.$ Since all
assumptions of the Condition 4.1 are satisfied, the above equation has a
unique solution 
\begin{equation*}
u\in W_{p,q}^{\left( 1,2\right) }(R\times \Omega )=\left\{ u\left(
t,x\right) \in X:\partial _{t}u,\partial _{xx}u\in X,\text{ and }%
u(t,x)_{x\in \partial \Omega }=0\right\} 
\end{equation*}
satisfying coercive estimate%
\begin{equation*}
\left\Vert \partial _{t}u\right\Vert _{X}+\left\Vert e^{-m\left\vert \cdot
\right\vert }\ast \partial _{t}u\right\Vert _{X}+\left\Vert \partial
_{xx}u\right\Vert _{X}+\left\Vert e^{-k\left\vert \cdot \right\vert }\ast
\partial _{xx}u\right\Vert _{X}\leq C\left\Vert f\right\Vert _{X}.
\end{equation*}

\section*{5. Remarks on Parabolic and Elliptic DOE}

In recent years Lutz Weis [26] \ and Herbert Amann [1] established maximal $%
L_{p}(E)$ and $B_{q,r}^{s}(E)$ regularity for abstract cauchy problem%
\begin{equation*}
\ 
\begin{array}{l}
\displaystyle u^{\prime }(t)+Au=f(t) \\ 
\\ 
\displaystyle u(0)=0%
\end{array}%
\eqno(5.1)
\end{equation*}%
Here based on obtained FMT, we shall give some remarks on $(5.1)$.

\vspace{3mm}

\textbf{Remark 5.1. }Let $E$ be an UMD space$.$ Suppose $A$ is an $R$%
-positive operator in $E$ i.e. 
\begin{equation*}
R\left\{ (1+\lambda )R(\lambda ,A):\text{ }\lambda \in S_{\varphi }\text{
for }\frac{\pi }{2}<\varphi <\pi \right\} <\infty .\eqno(5.2)
\end{equation*}%
Then for each $f\in L_{q}(R^{+};E),~(5.1)$ has a unique solution 
\begin{equation*}
u\in \bigcap\limits_{q<\theta <\infty }L_{\theta }(R^{+};E)\tbigcap
W_{q}^{1}(R^{+},E(A),E)
\end{equation*}%
satisfying coercive estimate%
\begin{equation*}
\left\Vert u\right\Vert _{W_{q}^{1}(R^{+},E(A),E)}+\left\Vert u\right\Vert
_{L_{\theta }(R^{+};E)}\leq C\left\Vert f\right\Vert _{L_{q}(R^{+};E)}.\eqno%
(5.3)
\end{equation*}%
Since $A$ is a generator of bounded analytic semigroup $T_{t}$, solutions of
(5.1) can be represented in the form of%
\begin{equation*}
\begin{array}{lll}
u(t) & = & \displaystyle\dint\limits_{0}^{t}T_{t-s}\text{ }f(s)\text{ }ds%
\end{array}%
\end{equation*}%
where%
\begin{equation*}
\begin{array}{lll}
T_{t-s} & = & \displaystyle T(t-s)=e^{-A(t-s)}.%
\end{array}%
\end{equation*}%
Therefore,%
\begin{equation*}
\begin{array}{lll}
u^{\prime }(t) & = & \displaystyle f(t)-\dint\limits_{0}^{t}AT_{t-s}\text{ }%
f(s)\text{ }ds.%
\end{array}%
\end{equation*}%
and%
\begin{equation*}
\begin{array}{lll}
Au(t) & = & \displaystyle\dint\limits_{0}^{t}AT_{t-s}\text{ }f(s)\text{ }ds.%
\end{array}%
\end{equation*}%
Now, it is easy to see that maximal $L_{q}(R^{+};E)$ to $L_{\theta
}(R^{+};E) $ regularity of (5.1) is equivalent to the boundedness of operator%
\begin{equation*}
Kf(t)=(AT)\ast f=\dint\limits_{-\infty }^{\infty }AT_{t-s}(f(s))ds=F^{-1} 
\left[ \left( AT\right) ^{\symbol{94}}(\cdot )\hat{f}(\cdot )\right]
\end{equation*}%
where 
\begin{equation*}
AT(t)=\left\{ 
\begin{array}{ll}
AT_{t}\text{ }\mbox{for} & t>0 \\ 
0\text{ }\mbox{for} & t\leq 0%
\end{array}%
\right. \text{.}
\end{equation*}

In order to show%
\begin{equation*}
\begin{array}{lll}
\left\Vert u\right\Vert _{L_{\theta }(R^{+};E)} & \leq & \displaystyle %
C\left\Vert f\right\Vert _{L_{q}(R^{+};E)},\text{ }%
\end{array}%
\end{equation*}%
and%
\begin{equation*}
\begin{array}{lll}
\left\Vert u\right\Vert _{W_{_{\theta }}^{1}(R^{+},E(A),E)} & \leq & %
\displaystyle C\left\Vert f\right\Vert _{L_{q}(R^{+};E)},\text{ for }1<q\leq
\theta <\infty%
\end{array}%
\end{equation*}%
it suffices to prove%
\begin{equation*}
m_{0}(t)=(T)^{\symbol{94}}(t)=R(it,A)
\end{equation*}%
and 
\begin{equation*}
m_{1}(t)=(AT)^{\symbol{94}}(t)=AR(it,A)\left( t\right) =itR(it,A)-I
\end{equation*}%
are Fourier multipliers$.$ Hence, we have to show%
\begin{equation*}
R\left\{ |t|^{\frac{1}{q}-\frac{1}{\theta }}m_{i}(t)\mid t\in R\backslash
\left\{ 0\right\} \right\} \leq C_{1}\eqno(5.4)
\end{equation*}%
and%
\begin{equation*}
R\left\{ |t|^{1+\frac{1}{q}-\frac{1}{\theta }}\frac{d}{dt}m_{i}(t)\mid t\in
R\backslash \left\{ 0\right\} \right\} \leq C_{2},\eqno(5.5)
\end{equation*}%
for $i=0,1.$ For the function $m_{0},$ (5.4) and (5.5) hold for each $\theta 
$ satisfying $1<q\leq \theta <\infty $ due to (5.2). However for the second
function we have 
\begin{equation*}
\left\Vert |t|^{\frac{1}{q}-\frac{1}{\theta }}m_{1}\right\Vert \leq C|t|^{1+%
\frac{1}{q}-\frac{1}{\theta }}\left\Vert R(it,A)\right\Vert \leq C\frac{%
|t|^{1+\frac{1}{q}-\frac{1}{\theta }}}{1+\left\vert t\right\vert }
\end{equation*}%
which implies that the right hand side is unbounded whenever $q<\theta .$
Thus (5.4) and (5.5) do not hold for $m_{1}$ unless $q=\theta $. Eventually, 
$(5.1)$ has a unique solution 
\begin{equation*}
u\in \bigcap\limits_{q<\theta <\infty }L_{\theta }(R^{+};E)\tbigcap
W_{q}^{1}(R^{+},E(A),E)
\end{equation*}%
satisfying the coercive estimate (5.3). \ \ \hbox{\vrule
height7pt width5pt}

\vspace{3mm}

\textbf{Remark 5.2.} Let $E$ be an UMD space and $A$ be an $R$-positive
operator in $E$. Then for each $f\in L_{q}(R;E),~$the following elliptic DOE%
\begin{equation*}
-u^{\prime \prime }(t)+Au=f(t)\eqno(5.6)
\end{equation*}%
has a unique solution 
\begin{equation*}
u\in \bigcap\limits_{q<\theta <\infty }W_{\theta }^{1}(R;E)\tbigcap
W_{q}^{2}(R;E(A),E)
\end{equation*}%
satisfying coercive estimate%
\begin{equation*}
\left\Vert u\right\Vert _{W_{q}^{2}(R;E(A),E)}+\left\Vert u\right\Vert
_{W_{\theta }^{1}(R;E)}\leq C\left\Vert f\right\Vert _{L_{q}(R;E)}.
\end{equation*}%
Applying Fourier transform to equation $(5.6)$, we obtain 
\begin{equation*}
\lbrack \xi ^{2}+A]\hat{u}(\xi )=\hat{f}(\xi ).
\end{equation*}%
Since $A$ is $R$-positive, solutions of $(5.6)$ can be represented in the
following form 
\begin{equation*}
u(x)=F^{-1}[A+\xi ^{2}]^{-1}\hat{f}.\eqno(5.7)
\end{equation*}%
By using $(5.7)$, we get 
\begin{equation*}
\begin{array}{lll}
\Vert u^{\prime }\Vert _{L_{\theta }(R;E)} & = & \displaystyle\left\Vert
F^{-1}\left[ \xi \left( A+\xi ^{2}\right) ^{-1}\hat{f}\right] \right\Vert
_{L_{\theta }(R;E)} \\ 
&  &  \\ 
\Vert Au\Vert _{L_{\theta }(R;E)} & = & \displaystyle\left\Vert F^{-1}\left[
A(A+\xi ^{2})^{-1}\hat{f}\right] \right\Vert _{L_{\theta }(R;E)} \\ 
&  &  \\ 
\Vert u^{\prime \prime }\Vert _{L_{\theta }(R;E)} & = & \displaystyle%
\left\Vert F^{-1}[\xi ^{2}(A+\xi ^{2})^{-1}\hat{f}]\right\Vert _{L_{\theta
}(R;E)}.%
\end{array}%
\end{equation*}%
Therefore, it suffices to show operator--functions 
\begin{equation*}
\sigma _{0}(\xi )=[A+\xi ^{2}]^{-1},\text{ }\sigma _{1}(\xi )=\xi (A+\xi
^{2})^{-1},\text{ }\sigma _{2}(\xi )=\xi ^{2}(A+\xi ^{2})^{-1}
\end{equation*}%
and 
\begin{equation*}
\sigma _{3}(\xi )=A[A+\xi ^{2}]^{-1}=I-\xi ^{2}[A+\xi ^{2}]^{-1}
\end{equation*}%
are Fourier multipliers$.$ Applying similar techniques as in the Remark 5.1
one can easily show that $\sigma _{j}$ are FMT. 
\hbox{\vrule height7pt
width5pt}\vspace{3mm}

\textbf{Example 5.3. }As an application of the Remark 5.1 we can give a
mixed problem for infinite system of diffusion equations i.e.%
\begin{equation*}
\begin{array}{l}
\displaystyle\frac{\partial u_{k}(t,x)}{\partial t}-(\Delta
+c)u_{k}(t,x)=f_{k}(t,x), \\ 
\\ 
\displaystyle u_{k}(0,x)=0,\text{ }u_{k}(t,x)|_{x\in \partial G}=0%
\end{array}%
\text{for }k=0,1\cdot \cdot \cdot .\eqno(5.8)
\end{equation*}%
We assume $c>0,$ $1<q<\infty ,$ $p\in \left( 1,\infty \right) $, $%
E=L_{p}(G;l_{p}),$ $A=-(\Delta +c)$ and $E(A)=W_{p}^{2}(G;l_{p})$. Here $%
G\subseteq R^{n}$ and $u\in W_{p}^{2}(G;l_{p})$ is assumed to satisfy a
boundary condition $u|_{\partial G}=0.$ Then, for each $f\in L_{q}(R^{+};E)$
(5.8) has a unique solution 
\begin{equation*}
u\in \bigcap\limits_{q<\theta <\infty }L_{\theta
}(R^{+};L_{p}(G;l_{p}))\tbigcap W_{q}^{1}(R^{+},E(A),E)
\end{equation*}%
and the following coercive estimate holds%
\begin{equation*}
\left\Vert u\right\Vert _{W_{q}^{1}(R^{+},E(A),E)}+\left\Vert u\right\Vert
_{L_{\theta }(R^{+};\text{ }L_{p}(G;l_{p}))}\leq C\left\Vert f\right\Vert
_{L_{q}(R^{+};\text{ }L_{p}(G;l_{p}))}.
\end{equation*}%
\textbf{Remark 5.4. }The FMT play an important role in the study of\
embedding theorems. For instance, under certain abstract conditions on a
Banach space $E,$ author proved in [20] that 
\begin{equation*}
D^{\alpha }:W_{p}^{l}\left( \Omega ;E\left( A\right) ,E\right) \rightarrow
L_{q}\left( \Omega ;E\left( A^{1-x}\right) \right)
\end{equation*}%
is continuous for $x=\frac{\alpha +\frac{1}{p}-\frac{1}{q}}{l}\leq 1,$ and
the Gagilardo-Nirenberg type sharp estimate 
\begin{equation*}
\left\Vert D^{\alpha }u\right\Vert _{L_{q}\left( \Omega ;E\left(
A^{1-\varkappa -\mu }\right) \right) }\leq h^{\mu }\left\Vert u\right\Vert
_{W_{p}^{l}\left( \Omega ;E\left( A\right) ,E\right) }+h^{-\left( 1-\mu
\right) }\left\Vert u\right\Vert _{L_{p}\left( \Omega ;E\right) },
\end{equation*}%
holds for all $u\in W_{p}^{l}\left( \Omega ;E\left( A\right) ,E\right) ,$ $%
1<p\leq q<\infty $, $0<\mu \leq 1-x$ and $0<h\leq h_{0}<\infty .$ However,
using the same techniques as in [20] and applying Theorem 2.3 one can remove
this combined assumption on $E$. Note that these embedding theorems play key
role in the theory of DOE, especially in estimation of lower order terms in
DOE of type e.g.%
\begin{equation*}
~-u^{\prime \prime }(t)+A_{1}(t)u^{\prime }(t)+Au(t)~=~f(t)
\end{equation*}%
where $A_{1}(t)$ is a variable and generally unbounded operator. More
general form of above equation and the embedding theorems are studied in
[21-22].

\begin{center}
\bigskip

\textbf{Acknowledgements}
\end{center}

\textit{The author would like to express a deep gratitude to Reviewer of his
previous submission who suggested him a way to prove Theorem 2.2.}

\end{document}